%
\magnification=\magstep1
\vsize=23.5truecm

\input amstex
\UseAMSsymbols
\input pictex
\NoBlackBoxes
  \font\gross=cmbx10 scaled\magstep1 
   \font\rmk=cmr8    \font\itk=cmti8  \font\ttk=cmtt8

   \newcount\notenumber
   
   \def\note{\advance\notenumber by 1 
       \plainfootnote{$^{\the\notenumber}$}}  
\def\mod{\operatorname{mod}}

\def\op{{\text{op}}}

\def\Tr{\operatorname{Tr}}

\def\Hom{\operatorname{Hom}}

\def\End{\operatorname{End}}
\def\Ext{\operatorname{Ext}}

\def\add{\operatorname{add}}
 
\def\tors{\operatorname{tors}} 
\def\ftors{\operatorname{f-tors}}

\def\Rahmen#1%
   {$$\vbox{\hrule\hbox%
                  {\vrule%
                       \hskip0.5cm%
                            \vbox{\vskip0.3cm\relax%
                               \hbox{$\displaystyle{#1}$}%
                                  \vskip0.3cm}%
                       \hskip0.5cm%
                  \vrule}%
           \hrule}$$}

   
	\vglue1truecm
\plainfootnote{}
{\rmk 2010 \itk Mathematics Subject Classification. \rmk 
Primary 
        16G10, 
        18E40. 
Secondary:
        05E10, 
        16D90, 
        16G70. 
}
\centerline{\gross Lattice structure of torsion classes for hereditary artin algebras.}     
		   		       \medskip\smallskip
\centerline{Claus Michael Ringel}
		  	    \bigskip\medskip

{\narrower\narrower
Abstract: Let $\Lambda$ be a connected hereditary artin algebra. We show that the set of functorially finite 
torsion classes of $\Lambda$-modules is a lattice if and only if $\Lambda$ is either representation-finite 
(thus a Dynkin algebra) or $\Lambda$ has only two simple modules. For the case of $\Lambda$ being the path algebra 
of a quiver, this result has recently been established by Iyama-Reiten-Thomas-Todorov and our proof follows closely
their considerations.\par} 
      \bigskip
Let $\Lambda$ be a connected hereditary artin algebra. The modules considered here
are left $\Lambda$-modules of finite length, $\mod\Lambda$ denotes the corresponding category.
The subcategories of $\mod\Lambda$ we deal with are always assumed to be closed under direct sums and direct
summands (in particular closed under isomorphisms). In this setting, a subcategory is a {\it torsion
class} (the class of torsion modules for what is called a torsion pair or a torsion theory)
provided it is closed under factor modules and extensions. The torsion classes form a partially ordered set
with respect to inclusion, it will be denoted by $\tors \Lambda$. This poset clearly is a lattice (even a
complete lattice). Auslander and Smal\o{} have pointed out that a torsion class $\Cal C$ in $\mod\Lambda$ 
is functorially finite
if and only if it has a cover (a {\it cover} for $\Cal C$ 
is a module $C$ such that $\Cal C$ is the set of modules generated by $C$),
we denote by $\ftors \Lambda$ the set of functorially finite torsion classes in $\mod\Lambda$.

In a recent paper [IRTT], Iyama, Reiten, Thomas and Todorov have discussed the question whether also the poset
$\ftors \Lambda$ (with the inclusion order) is a lattice. 
	\medskip 
{\bf Theorem.} {\it The poset $\ftors\Lambda$ is a lattice if and only if $\Lambda$ is representation finite or
$\Lambda$ has precisely two simple modules.}
	\medskip 
Iyama, Reiten, Thomas, Todorov have shown this in the special case when $\Lambda$ is a $k$-algebra with
$k$ an algebraically closed field (so that $\Lambda$ is Morita equivalent to the path algebra of a quiver).
The aim of this note is to provide a proof in general. We follow closely the strategy of the paper [IRTT]
and we will use Remark 1.13 of [IRTT] which asserts
that a meet or a join of two elements $\Cal C_1, \Cal C_2$
in $\ftors \Lambda$ exists if and only if the meet or the join of $\Cal C_1, \Cal C_2$ formed in $\tors\Lambda$
belongs to $\ftors\Lambda$, respectively. 
	\bigskip\medskip 
{\bf 1. Normalization.}
	\medskip 
Let $\Cal X$ be a class of modules. We denote by $\add(\Cal X)$ the modules
which are direct summands of direct sums of modules in $\Cal X$. 
A module $M$ is {\it generated} by $\Cal X$ provided $M$ is a factor module of a module in $\add(\Cal X)$,
and $M$ is {\it cogenerated} by $\Cal X$ provided $M$ is a submodule of a module in $\add(\Cal X)$. The 
subcategory of all modules generated by $\Cal X$ is denoted by $\Cal G(\Cal X).$ In case $\Cal X = \{X\}$
or $\Cal X = \add X$, we write $\Cal G(X)$ instead of $\Cal G(\Cal X)$, and use the same
convention in similar situations. We write $\Cal T(X)$ for the smallest torsion class containing the module $X$
(it is the intersection of all torsion classes containing $X$, and it can be constructed as the closure of $\{X\}$
using factor modules and extensions). 

Since $\Lambda$ is assumed to be hereditary, we write $\Ext(X,Y)$ instead of $\Ext^1(X,Y).$ 
Recall that a module $X$ is said to be {\it exceptional} provided it is indecomposable and has no
self-extensions (this means that $\Ext(X,X) = 0$).
	\medskip 

Following Roiter [Ro], we say that a module $M$ is {\it normal} provided there is no proper direct decomposition
$M = M'\oplus M''$ such that $M'$ generates $M''$ (this means: if $M = M'\oplus M''$ and $M'$ generates $M''$, then
$M'' = 0$). Of course, given a module $M$, there is a direct decomposition $M = M'\oplus M''$ such that $M'$ is
normal and $M'$ generates $M''$ and one can show that $M'$ is determined by $M$ uniquely up to isomorphism,
thus we call $M' = \nu(M)$ a {\it normalization} of $M$. This was shown already by Roiter [Ro], and later 
by Auslander-Smal\o{} [AS]. It is also a consequence of the following Lemma which will be needed for our
further considerations.
	\medskip
{\bf Lemma 1.} (a) {\it Let $(f_1,\dots,f_t,g)\:X \to X^t\oplus Y$ be an injective map for some natural number $t$,
with all the maps $f_i$ in the radical of $\End(X)$. Then $X$ is cogenerated by $Y$.}

(b) {\it Let $(f_1,\dots,f_t,g)\:X^t \oplus Y \to X$ be a surjective map for some natural number $t$,
with all the maps $f_i$ in the radical of $\End(X)$, then $Y$ generates $X$.}
	\medskip
Proof. (a) Assume that the radical $J$ of $\End(X)$ satisfies $J^m = 0$.
Let $W$ be the set of all compositions $w$ of at most $m-1$ maps of the form $f_i$ with $1\le i \le t$
(including $w = 1_X$).  
We claim that $(gw)_{w\in W}\:X \to Y^{|W|}$ is injective. Take a non-zero element $x$
in $X$. Then there is $w\in W$ such that $w(x)\neq 0$ and $f_iw(x) = 0$ for $1\le i\le t$.
Since $(f_1,\dots,f_t,g)$ in injective and $w(x) \neq 0$, we have $(f_1,\dots,f_t,g)(w(x)) \neq 0.$ 
But $f_iw(x) = 0$ for $1\le i \le t,$ thus $g(w(x)) \neq 0.$ This completes the proof.
	
(b) This follows by duality.    \hfill$\square$
	\bigskip 

{\bf Corollary (Uniqueness of normalization).} 
{\it Let $M$ be a module. Assume that $M = M_0\oplus M_1 = M_0'\oplus M_1'$ such that both 
$M_0$ and $M'_0$ generate $M$. Then there is a module $N$ which is a direct summand of both $M_0$ and
$M'_0$ which generates $M$.}
	\medskip
Proof: We may assume that $M$ is multiplicity free. Write 
$M_0 \simeq N\oplus C,\ M_0' \simeq N \oplus C',$ such that $C, C'$
have no indecomposable direct summand in common. 
Now, $N\oplus C$ generates $N\oplus C'$ generates $N\oplus C$ generates $C$.
We see that $N\oplus C$ generates $C$, such that the maps $C \to C$ used 
belong to the radical of $\End(C)$ (since they factor through $\add(N\oplus C')$
and no indecomposable direct summand of $C$ belongs to $\add(N\oplus C')$).
Lemma 1 asserts that $N$ generates $C$, thus it generates $M$.
  \hfill$\square$

	\bigskip
{\bf Proposition 1.} {\it If $T$ has no self-extensions, then $T$ is a cover  
for the torsion class $\Cal T(T)$.
Conversely, if $\Cal T$ is a torsion class with cover $C$, then $\nu(C)$ has no
self-extensions.}
	\medskip
Proof. For the first assertion, one has to observe that $\Cal G(T)$ is closed under extensions, 
thus equal to $\Cal T(T).$ This is a standard result say in tilting theory. Here is the argument:
let $g'\:T' \to M'$ and $g''\:T'' \to M''$ be surjective maps with $T',T''$ in $\add T$. 
Let $0\to M' \to M \to M'' \to 0$ be an exact sequence. The induced exact sequence with respect to $g''$
is of the form $0 \to M' \to Y_1 \to T'' \to 0$ with a surjective map $g_1\:Y_1 \to M$. 
Since $\Lambda$ is hereditary and $g'$ is surjective, there is an exact sequence
$0 \to T' \to Y_2 \to T'' \to 0$ with a surjective map $g_2\:Y_2 \to Y_1.$ Since $\Ext(T'',T') = 0$,
we see that $Y_2$ is isomorphic to $T'\oplus T''$, thus in $\add T$. And there is the surjective map 
$g_1g_2\:Y_2 \to M.$
	
For the converse, we may assume that $C$ is normal and have to show that $C$ has no self-extension. 
Let $C_1, C_2$ be indecomposable direct summands of $C$ and assume for the contrary that there is 
a non-split exact sequence
$$
 0 \to C_1 \to M \to C_2 \to 0.
$$
Now $M$ belongs to $\Cal T$, thus it is generated by $C$, say 
there is a surjective map $C' \to M$ with $C'\in \add C.$ Write $C' = C_2^t\oplus C''$
such that $C_2$ is not a direct summand of $C''$. Consider the surjective map
$C_2^t\oplus C'' \to M \to C_2$. Since the last map $M \to C_2$ is not a split epimorphism,
all the maps $C_2 \to C_2$ involved  belong to the radical of $\End(C_2)$.
According to Lemma 1, $C''$ generates $C_2.$ This contradicts the assumption that $C$
is normal.    \hfill$\square$
	\medskip 
{\bf Remark.} As we have mentioned, normal modules have been considered by Roiter, but actually, 
he used a slightly deviating name, calling them  ''normally indecomposable''.
	\bigskip\medskip 
{\bf 2. $\Ext$-cycles.}
	\medskip 
An {\it $\Ext$-cycle} of cardinality $t$ is a sequence 
$X_1,X_2,\dots, X_m$ of pairwise orthogonal bricks
such that $\Ext(X_{i-1}, X_i) \neq 0$ for $1\le i \le m$, 
with $X_0 = X_m$.
An {\it $\Ext$-pair} is  an $\Ext$-cycle of cardinality $2$
consisting of exceptional modules. (One may call an $\Ext$-cycle $X_1,\dots, X_m$ 
{\it minimal} provided there is no $\Ext$-cycle of smaller cardinality which uses
(some of) these modules. Using this definition, the $\Ext$-pairs are just the
minimal $\Ext$-cycles of cardinality $2$.)
	\bigskip 
{\bf Proposition 2.} {\it If $X_1,X_2,\dots, X_m$ is an
$\Ext$-cycle, then $\Cal T(X_1,\dots,X_m)$ has no cover.}
	\medskip 
Proof: Let $\Cal F =
\Cal F(X_1,\dots,X_m)$ be the extension
closure of $X_1,\dots,X_m$, 
thus the class of modules with a filtration with factors of the form $X_i$, where $1\le i \le m$. 
According to [R], $\Cal F$ is
an abelian subcategory with exact embedding functor, with (relative) simple objects 
the modules $X_1,\dots,X_m$. The objects in $\Cal F$ have finite (relative) length, thus
also the (relative) Loewy length for these objects is defined. 
We denote by $\Cal F_t$ the full subcategory of objects in $\Cal F$ of (relative) Loewy
length at most $t$. 

We have 
$$
 \Cal F_1 \subseteq \Cal F_2 \subseteq \cdots \subseteq \Cal F_t \subseteq \cdots
$$
and therefore 
$$
 \Cal G(\Cal F_1) \subseteq \Cal G(\Cal F_2) \subseteq \cdots \subseteq \Cal G(\Cal F_t) 
 \subseteq \cdots,
$$
Let $\Cal G = \bigcup\nolimits_t \Cal G(\Cal F_t)$. We claim that 
$\Cal G =  \Cal T(X_1,\dots,X_m).$ The modules in $\Cal G$ belong to 
$\Cal T(X_1,\dots,X_m)$ and $X_1,\dots, X_m$ belong to $\Cal G.$ 
Thus, it is sufficient to show that $\Cal G$ is a torsion class.

Since $\Cal G$ is the filtered union of classes
closed under epimorphisms, it is closed under epimorphisms. 
In order to show that $\Cal G$ is closed under extensions, we follow the proof 
for the first assertion of Proposition 1 as closely as possible:
Let $g'\:F' \to M'$ and $g''\:F'' \to M''$ be surjective maps with $F',F''$ in $\add \Cal F_s$ for some $s$. 
Let $0\to M' \to M \to M'' \to 0$ be an exact sequence. The induced exact sequence with respect to $g''$
is of the form $0 \to M' \to Y_1 \to F'' \to 0$ with a surjective map $g_1\:Y_1 \to M$. 
Since $\Lambda$ is hereditary and $g'$ is surjective, there is an exact sequence
$0 \to F' \to Y_2 \to F'' \to 0$ with a surjective map $g_2\:Y_2 \to Y_1.$ Since $F', F''$ belong to
$\Cal F$ and their (relative) Loewy length is at most $s$, 
the exact sequence shows that $M$ also belongs to $\Cal F$ and has (relative)
Loewy length at most $2s$.
The surjective map $g_1g_2\:Y_2 \to M$ shows that $M$ is in $\Cal G(\Cal F_{2s}) \subseteq \Cal G$.
	
Now assume that $C$ is a cover for $\Cal G$. The module $C$ belongs to $\Cal G(\Cal F_r)$ for some $r$,
thus there is an epimorphism $f\:F \to C$ for some $F \in \Cal F_r$. With $C$ also $F$
is a cover for $\Cal G$. Note that there is a module $F'$ which belongs to $\Cal F_{r+1}$ 
and not to $\Cal F_r$, for example any object in $\Cal F$ which is 
(relative) serial  and has (relative) length equal to $r+1$.
Since $F'$ is in $\Cal G$, and $F$ is a cover of $\Cal G$, the module $F'$ is generated by $F$. 
But if $F'$ is generated by $F$, its (relative) Loewy length is at most $r$. 
This means that $F'$ is in $\Cal F_r$, a contradiction. 
	\hfill $\square$
	\bigskip\medskip 
{\bf 3. Construction of $\Ext$-pairs.}
	\medskip 
{\bf Proposition 3.} {\it A connected hereditary artin algebra which is representation-infinite and has 
at least three simple modules has $\Ext$-pairs.}
	\medskip 
Given a finite dimensional algebra $R$, we denote by $Q(R)$ its {\it $\Ext$-quiver:} its
vertices are the isomorphism classes $[S]$ of the simple $R$-modules  $S$, 
and given two simple $R$-modules $S,S'$, there
is an arrow $[S] \to [S']$ provided $\Ext^1(S,S') \neq 0.$ If $R$ is hereditary, then 
clearly $Q(R)$ is directed. If necessary, we endow $Q(R)$ with a valuation as 
follows: Given an arrow $S \to S'$, consider $\Ext(S,S')$ as a left $\End(S)^\op$-module
or as a left $\End(S')$-module and put
$$
 v([S],[S']) = (\dim {}_{\End(S)}\Ext(S,S'))(\dim {}_{\End(S')^\op}\Ext(S,S'))
$$
(note that in contrast to [DR], we only will need the product of the two dimensions, not the pair).  
Given a vertex $i$ of $Q(R)$, we denote by $S(i), P(i), I(i)$
a simple, projective or injective module corresponding to the vertex $i$, respectively.

We later will use the following: 
If $Q(\Lambda) = (1\to 2)$, then the arrow $1\to 2$ has valuation at least $2$ if and
only if $I(2)$ is not projective if and only if $P(1)$ is not injective;
if the arrow $1\to 2$ has valuation at least $3$, then $\tau S(1)$ (where $\tau$ is the Auslander-Reiten
translation) is neither projective, nor a neighbor of $P(1)$ in the Auslander-Reiten quiver,
consequently $\Hom(P(1),\tau^2S(1)) \neq 0,$ thus $\Ext(\tau S(1),P(1))\neq 0.$ 

For any hereditary algebra $\Lambda$ with $Q(\Lambda)$ being a tree quiver, it is easy
to construct a sincere exceptional module, using induction: If $Q'$ is a subquiver
of $Q$ such that $Q$ is obtained from $Q'$ by adding just one vertex $\omega$ 
and one arrow, and
$M'$ is an exceptional module for the restriction of $\Lambda$ to $Q'$, then
let $M$ be the universal extension of $M'$ by copies of $S(\omega)$; here we consider
extensions from above or from below, provided $\omega$ is a source or a sink,
respectively.
	\medskip 
For the proof of Proposition 3, we consider four special cases:
	\medskip
{\bf Case 1. The algebra $\Lambda$ is tame.}
	
We use the structure of the Auslander-Reiten quiver of $\Lambda$ as presented in [DR].
Since we assume that $\Lambda$ has at least 3 vertices, there is a tube of rank $r \ge 2$.
The simple regular modules in this component
form an $\Ext$-cycle of cardinality $r$, say $X_1,\dots, X_r.$ 
There is a unique indecomposable module $Y$ with a filtration 
$Y = Y_0 \supset Y_1 \supset \cdots \supset Y_{r-1} = 0$ such that 
$Y_{i-1}/Y_i = X_i$ for $1\le i \le r-1.$ 
Clearly, the pair $Y,X_r$ is an $\Ext$-pair.  

	\medskip 
{\bf Case 2. The quiver $Q = Q(\Lambda)$ is not a tree.}

Deleting, if necessary, vertices, we may assume that the underlying graph of $Q$ is
a cycle. 
Let $w$ be a path from a sink $i$ to a source $j$ of smallest length, let $Q'$ be the
subquiver of $Q$ given by the vertices and the arrows which occur in $w$. 
Not every vertex of $Q$ belongs to $Q'$, since otherwise $Q$ is obtained from $Q'$
by adding just arrows, thus by adding a unique arrow, namely an arrow $i\to j$.
But then this arrow is also a path from a sink to a source, and it has length $1$.
By the minimality of $w$, we see that also $w$ has length $1$ 
and therefore $Q$ has just the two vertices $i,j$. But then $Q$ can have only one
arrow, thus is a tree. This is a contradiction.

Let $Q''$ be the full subquiver given by all vertices of $Q$ which do not belong to 
$Q'$. Of course, $Q''$ is connected (it is a quiver of type $\Bbb A$). 
Let $X$ be an exceptional module with support $Q'$ and $Y$ an exceptional
module with support $Q''$. Since $Q', Q''$ have no vertex in common, we see that $\Hom(X,Y) = 0 
= \Hom(Y,X)$.

There is an arrow $i\to j''$ with $j''$ a vertex of $Q''$. This
arrow shows that $\Ext^1(X,Y) \neq 0.$ Similarly, 
there is an arrow $i''\to j$ with $i''$ a vertex of $Q''$. This
arrow shows that $\Ext^1(Y,X) \neq 0.$ 
	\bigskip 
We consider now algebras $\Lambda$ with $\Ext$-quiver $1\to 2 \to 3$. 
We denote by $\Lambda'$ the restriction of $\Lambda$ to the subquiver with vertices $1,2$,
and by $\Lambda''$ the restriction of $\Lambda$ to the subquiver with vertices $2,3.$
Given a representation $M$, let $M_3$ be the
sum of all submodules of $M$ which are isomorphic to $S(3),$ then $M/M_3$ is a $\Lambda'$-module. 
	\medskip 
{\bf Lemma 2.} {\it Let $X, Y$ be $\Lambda$-modules. If $X_3 = 0$ and 
$\Ext^1(Y/Y_3,X) \neq 0$, then also $\Ext^1(Y,X) \neq 0.$}
	\medskip
Proof: The exact sequence $0 \to Y_3 \to Y \to Y/Y_3 \to 0$ 
yields an exact sequence
$$
  \Hom(Y_3,X) \to \Ext^1(Y/Y_3,X)  \to \Ext^1(Y,X) 
$$
The first term is zero, since $Y_3$ is a sum of copies of $S(3)$ and $X_3 = 0$.
Thus, the map $\Ext^1(Y/Y_3,X)  \to \Ext^1(Y,X) $ is injective.
	\medskip
{\bf Case 3. $Q(\Lambda) = (1 \to 2 \to 3)$, and $v(1,2)\ge 2,\ v(2,3)\ge 2.$}
	\medskip
Let $X = S(2)$ and let $Y$ be the universal extension of $X$ using the modules $(1)$ and $S(3)$
(thus, we form the universal extension from above using copies of $S(1)$ and the universal
extension from below using copies of $S(3)$. Clearly, $Y$ is exceptional. Since the socle of $Y$
consists of copies of $S(3)$, we have $\Hom(S(2),Y) = 0.$ Since the top of $Y$
consists of copies of $S(1)$, we have $\Hom(Y,S(2)) = 0.$

Since $v(1,2) \ge 2$, the module $Y/Y_3$ is not a projective $\Lambda'$-module.
As a consequence, $\Ext(Y/Y_3,S(2)) \neq 0.$ Lemma 2 shows that
also $\Ext(Y,S(2)) \neq 0.$ 
By duality, we similarly see that $\Ext(S(2),Y) \neq 0.$
	\medskip
{\bf Case 4. $Q(\Lambda) = (1 \to 2 \to 3)$, and $v(1,2)\ge 3,\ v(2,3) = 1.$}
	\medskip
Let $X = P(1)/P(1)_3$ (thus $X$ is the projective $\Lambda'$-module with top $S(1)$).
Let $Y = \tau X$, where $\tau = D\Tr$ is the Auslander-Reiten translation in $\mod\Lambda$.
Of course, both
modules $X,Y$ are exceptional. Since $Y = \tau X,$ we know already that $\Ext^1(X,Y) \neq 0.$

We claim that $Y/Y_3 = \tau'S(1)$, where $\tau'$ is the Auslander-Reiten translation of $\Lambda'$.
Since $P(1)_3 = S(3)^a$ for some $a \ge 1$, a minimal projective presentation of $X$ has the form 
$$
 0 \to S(3)^a \to P(1) \to X \to 0, \tag{$*$}
$$
thus the defining exact sequences for $Y = \tau X$ is of the form
$$
 0 \to Y \to I(3)^a \to S(1) \to 0.
$$
In order to obtain $\tau'S(1)$, we start with a minimal projective presentation
$$
 0 \to S(2)^a \to P'(1) \to S(1) \to 0, \tag{$**$}
$$
where $P'(1)$ is the projective cover of $S(1)$ as a $\Lambda'$-module (actually, $P'(1)= X$). Since $\nu(2,3) = 1$, 
the number $a$ in $(*)$ and $(**)$ is the same. 
The defining exact sequences for $Y = \tau X$ and
$\tau'S(1)$ are part of the following commutative diagram with exact rows and columns:
$$
\CD
 @.      0        @.    0  \cr
 @.      @AAA          @AAA        \cr
 0 @>>> \tau'S(1) @>>> I(2)^a @>>> S(1) @>>> 0 \cr
 @.      @AAA          @AAA        @| \cr
 0 @>>>  Y        @>>> I(3)^a @>>> S(1) @>>> 0 \cr
 @.      @AAA          @AAA        \cr
 @.      S(3)^a  @=     S(3)^a \cr
 @.      @AAA          @AAA        \cr
 @.      0        @.    0  
\endCD
$$
The left column shows that $Y/Y_3 = \tau'S(1)$. 

We have noted already that $v(1,2) \ge 3$ implies that $\Ext(\tau' S(1),P'(1))\neq 0.$ 
According to Lemma 2, we see that $\Ext(Y,X) \neq 0.$ 

Finally, let us show that $X,Y$ are orthogonal. Any homomorphism $Y \to X$ vanishes on $Y_3$, since $X$ has no
composition factor $S(3)$. Now $Y/Y_3$ is indecomposable and not projective as a $\Lambda'$-module,
whereas $X$ is a projective $\Lambda'$-module, thus $\Hom(Y,X) = \Hom(Y/Y_3,X) = 0.$

On the other hand, the restriction $X''$ of $X$ to the subquiver $Q''$ with vertices $2,3$ is a sum of copies of $S(2)$,
whereas the restriction of $Y$ to the subquiver $Q''$ is a projective-injective module. It follows that
the restriction of any homomorphism $f\:X \to Y$ vanishes on $X''$. Thus $f$ factors through a direct sum
of copies of $S(1)$. But $S(1)$ is injective and obviously not a submodule of $Y$. It follows that $f = 0.$
	\medskip
{\bf Remark.} Concerning the cases 3 and 4, there is an alternative proof which uses dimension vectors and the 
Euler form on the Grothendieck group $K_0(\Lambda)$. But for this approach, one needs to deal with the
valuation of $Q(\Lambda)$ as in [DR], attaching to any arrow $i\to j$ a pair $(a,b)$ of positive numbers 
instead of the single number $v(i,j) = ab$. 
	\bigskip
Proof of Proposition 3. Let $\Lambda$ be connected, hereditary, representation-infinite, with at least 3 simple
modules. Case 2 shows that we can assume that $Q(\Lambda)$ is a tree. Assume that there is a subquiver $Q'$
such that at least two of the arrows have valuation at least $2$, choose such a $Q'$ of minimal length. 
We want to construct an $\Ext$-pair for the restriction of $\Lambda$ to $Q'$. Using reflection functors (see [DR]),
we can assume that $Q'$ has orientation $1 \to 2 \to \cdots \to n\!-\!1 \to n$. If $n = 3,$ then this is
case 3. Thus assume $n\ge 4.$ The minimality of $Q'$ asserts that $\nu(i,i+1) = 1$ for $2\le i \le n-2.$ 
If we denote by $\Lambda'$ the restriction of $\Lambda$ to $Q'$, then $\Lambda'$ has a full exact abelian
subcategory $\Cal U$ which is equivalent to the module category of an algebra as discussed in case 3 (namely the
subcategory of all $\Lambda'$-modules which do not have submodules of the form $S(i)$ with $2\le i \le n-2$
and no factor modules of the form $S(i)$ with $3\le i \le n-1$). Since $\Cal U$ has $\Ext$-pairs, also $\mod\Lambda$
has $\Ext$-pairs. Thus, we can assume that at most one arrow $i\to j$ has valuation greater than $2$. If $v(i,j) \ge 3$,
then we take a connected subquiver $Q'$ with 3 vertices containing this arrow $i\to j$. If necessary, we use again
reflection functors in order to change the orientation so that we are in case 4.
Thus we are left with the representation-infinite algebras $\Lambda$ with the following properties: 
$Q(\Lambda)$ is a tree, there is no arrow with valuation greater than $2$
and at most one arrow with valuation equal to $2$. It is easy to see that 
$Q(\Lambda)$ contains a subquiver $Q'$ such that the restriction of $\Lambda$ to $Q'$ is tame, thus we
can use case 1. \hfill$\square$
	\bigskip
Proof of Theorem. Let $\Lambda$ be connected and hereditary.
If $\Lambda$ is representation-finite, then $\tors \Lambda = \ftors\Lambda$, thus $\ftors\Lambda$ is a lattice.
If $\Lambda$ has precisely two simple modules, then $\ftors\Lambda$ can be described easily (see the proof of
Proposition 2.2 in [IRTT] which works in general), it obviously is a lattice.

On the other hand, if $\Lambda$ is representation-infinite and has at least three simple modules, then
Proposition 3 asserts that $\Lambda$ has an $\Ext$-pair, say $X,Y$. Since $X, Y$ are exceptional modules,
Proposition 1 shows that $\Cal T(X) = \Cal G(X)$ and $\Cal T(Y) = \Cal G(Y)$ both belong to $\ftors \Lambda$. 
The join of $\Cal T(X)$ and $\Cal T(Y)$ in $\tors \Lambda$ is $\Cal T(X,Y)$. According to Proposition 2, 
$\Cal T(X,Y)$ does not belong to $\ftors \Lambda$.   
\hfill$\square$

	\bigskip\medskip

{\bf 4. References}
	\medskip 
\item{[AS]} M.~Auslander, S.~O.~Smal\o: Preprojective modules over artin algebras,
    J. Algebra 66 (1980) 61--122.
\item{[DR]} V.~Dlab, C.~M.~Ringel:  Indecomposable representations of graphs and algebras. 
    Mem. Amer. Math. Soc. 173 (1976). 
\item{[IRTT]} O.~Iyama, I.~Reiten, H.~Thomas, G.~Todorov: Lattice structure of torsion classes for path algebras
   of quivers. arXiv:1312.3659
\item{[R]} C.~M.~Ringel: Representations of $k$-species and bimodules.  J. Algebra 41 (1976), 269--302. 
\item{[Ro]} A.~V.~Roiter: Unboundedness of the dimension of the indecomposable representations
   of an algebra which has infinitely many indecomposable representations. Izv.~Akad.~Nauk SSSR. Ser.~Mat.~32
   (1968), 1275-1282

	\bigskip\bigskip
{\rmk
C. M. Ringel\par
Department of Mathematics, Shanghai Jiao Tong University \par
Shanghai 200240, P. R. China, and \par 
King Abdulaziz University, P O Box 80200\par
Jeddah, Saudi Arabia\par


e-mail: {\ttk ringel\@math.uni-bielefeld.de} \par
}

\bye